# Main Manuscript for

## Global optimization tailored for graphics processing units: Complete and rigorous search for large-scale nonlinear minimization

Guanglu Zhang[a]*, Qihang Shan[a], Jonathan Cagan[a]

[a]Department of Mechanical Engineering, Carnegie Mellon University, Pittsburgh, PA 15213

*Corresponding author: Guanglu Zhang

**Email:** glzhang@cmu.edu

**Author Contributions:** G.Z. and J.C. designed research; G.Z., Q.S., and J.C. performed research; G.Z. and Q.S. analyzed data; G.Z. and J.C. wrote the paper; and Q.S. and G.Z. wrote code.

**Competing Interest Statement:** The authors declare no competing interest.

**Classification:** Physical Sciences/Computer Sciences

**This PDF file includes:**

Main Text
Figures 1 to 3
Tables 1 to 3

## Abstract

This paper introduces a numerical method to enclose the global minimum of a nonlinear function subject to simple bounds on the variables. Using interval analysis, coupled with the computational power and architecture of graphics processing units (GPUs), the method iteratively rules out the regions in the search domain where the global minimum cannot exist and leaves a finite set of regions where the global minimum must exist. For effectiveness, because of the rigor of interval analysis, the method is guaranteed to enclose the global minimum even in the presence of rounding errors. For efficiency, the method employs a novel GPU-based single program, single data parallel programming style to circumvent major GPU performance bottlenecks, and a variable cycling technique is also integrated into the method to reduce computational cost when minimizing large-scale nonlinear functions. The method is validated by minimizing 11 benchmark test functions with scalable dimensions, including the well-known Ackley function, Griewank function, Levy function, Rastrigin function, and Rosenbrock function. These benchmark test functions represent grand challenges of global optimization, and enclosing the guaranteed global minimum of these benchmark test functions with more than 80 dimensions has not been reported in the literature. Our method completely searches the feasible domain and successfully encloses the guaranteed global minimum of these 11 benchmark test functions with up to 10,000 dimensions using only one GPU in a reasonable computation time, far exceeding the reported results in the literature due to the unique method design and implementation based on GPU architecture.

## Significance Statement

A novel single program, single data (SPSD) parallel programming style is introduced for graphics processing unit (GPU) computing that circumvents two major GPU performance bottlenecks, CPU-GPU data transfers and GPU global memory access. By synthesizing interval analysis with the SPSD parallel programming style, the GPU-based numerical method presented in this paper can enclose the global minimum of large-scale multimodal nonlinear functions with more than 100 variables. Eleven well-known benchmark test functions with up to 10,000 variables are successfully minimized using the method. The GPU-based numerical method for large-scale nonlinear minimization serves as a critical tool that enables new scientific discoveries in many science and engineering fields.

## Main Text

## 1. Introduction

A nonlinear minimization problem is defined as

$$\underset{\boldsymbol{x}\in\mathbb{R}^n}{\text{minimize}} \quad f(\boldsymbol{x}), \tag{1}$$

$$\text{subject to} \quad \boldsymbol{l} \le \boldsymbol{x} \le \boldsymbol{u}, \tag{2}$$

where $\boldsymbol{l} \in \{ \mathbb{R} \cup \{-\infty\} \}^n$, $\boldsymbol{u} \in \{ \mathbb{R} \cup \{+\infty\} \}^n$, $\boldsymbol{l} < \boldsymbol{u}$, and $f : \mathbb{R}^n \to \mathbb{R}^1$. The nonlinear function $f$, also known as the objective function or loss function, has a closed-form expression, but it may or may not be continuous and differentiable. In addition, the nonlinear function is possibly highly nonconvex and may have many local minima. Of note, a nonlinear maximization problem could be transformed into the nonlinear minimization problem defined by Eq. (1) and Eq. (2) by adding a minus sign to the nonlinear function. Finding the global extremum (minimum or maximum) of a nonlinear function is also known as global optimization.

The nonlinear minimization problem defined by Eq. (1) and Eq. (2) is pervasive in many science and engineering fields. In many practical applications, it is desirable, and sometimes indispensable, to find the global minimum of a nonlinear function (1-3). However, popular numerical optimization methods, especially gradient-based methods (e.g., gradient descent methods (4, 5), interior-point methods (6, 7), trust region methods (8, 9), and BFGS methods (10, 11)), heuristic methods (e.g., genetic algorithms (12), simulated annealing (13), particle swarm optimization (14), and tabu search (15, 16)), and global optimization methods based on regular floating-point arithmetic (e.g., Nelder-Mead methods (17, 18), DIRECT methods (19, 20), basin-hopping (21), and differential evolution (22)) often fail to find the global minimum of a nonlinear function, because these popular numerical optimization methods suffer from at least one of the following three limitations:

- The convergence of some popular numerical optimization methods relies on an appropriate initial guess. The users of these methods must guess or randomly select the initial value(s) of each variable, also known as the starting point(s), to initialize these methods. A poor initial guess can contribute to the non-convergence of these methods or lead these methods to converge to the inferior solution (23, 24).
- These popular numerical optimization methods often become trapped in a local minimum. When the initial guess and the stopping criteria are properly specified, these methods may converge to a local minimum of the nonlinear function. However, the nonlinear function may have many local minima. Even if these methods are run multiple times with different initial guesses, the users are still not able to know whether the global minimum is found, and thus, it is difficult for the users to determine when to stop the search process using these methods.
- Many popular numerical optimization methods do not take rounding errors into account. When the nonlinear function is evaluated during the iteration process of these methods, the numerical result of each floating-point arithmetic operation is usually rounded to its closest binary machine number (25). Since the evaluation of a nonlinear function often consists of multiple floating-point arithmetic operations, the accumulation of rounding errors could lead to incorrect function evaluation results and wrong solutions (26, 27).

Based on the limitations of these popular numerical optimization methods, several global optimization methods using interval analysis (1, 26-30), referred to here as existing interval methods, have been introduced to enclose the global minimum of a nonlinear function. However, since these existing interval methods employ interval Newton methods and/or consistency techniques (e.g., box consistency and hull consistency) during their iteration process, they usually require the nonlinear function to be continuous and differentiable, and the users often need to transform the nonlinear function into a different form on a case-by-case basis. In addition, these existing interval methods are significantly more computationally expensive than the popular numerical optimization methods discussed above. As a result, these existing interval methods are not commonly adopted for minimizing nonlinear functions in practical applications, especially for minimizing large-scale nonconvex nonlinear functions that have more than 100 dimensions (i.e., the vector $x$ in Eq. (1) includes more than 100 variables) and many local minima. Importantly, the recent advances in graphics processing units (GPUs) have brought unprecedented computational power to users (31, 32), offering new opportunities to overcome the computational bottleneck in global optimization using interval analysis. Unfortunately, existing interval methods cannot utilize the computational power of GPUs, because modern GPUs are designed for massively parallel computing with a unique compute and memory architecture, but existing interval methods are designed for central processing unit (CPU)-based sequential computing, making these methods difficult to achieve massively parallel computing based on GPU architecture.

In this paper, we leverage interval analysis, coupled with the computational power and architecture of GPU, to design and implement a numerical method that is guaranteed to enclose the global minimum (or maximum) of a nonlinear function, especially a large-scale nonconvex nonlinear function, in a reasonable computation time. Instead of extending or accelerating any existing CPU-based numerical optimization methods using GPU, our new method is originally tailored for GPU, and the design and implementation of the method are intertwined based on GPU architecture. In other words, the method is a GPU-based method rather than a GPU-accelerated method. In contrast to the limitations of existing numerical optimization methods for minimizing nonlinear functions, our method is a general, complete, rigorous, and efficient global optimization method:

- *Generality*: The method does not require the nonlinear function in Eq. (1) to be continuous or differentiable, and the users do not need to transform the nonlinear function into a different form. The allowable range of each variable in Eq. (2) could be [-∞, +∞] because the method uses the closed interval system (27).
- *Completeness*: The method encloses the global minimum of a nonlinear function with certainty. The method finds all the *n*-dimensional regions (based on user-specified tolerances) in the search domain where the global minimum of the nonlinear function must exist after a finite computation time.
- *Rigorousness*: The method is guaranteed to enclose the global minimum of a nonlinear function even in the presence of rounding errors. The method also provides computer-aided proof for the existence of the global minimum of a nonlinear function in the search domain.
- *Efficiency*: The method is originally designed and implemented for GPU based on GPU compute and memory architecture. Importantly, the method employs a novel GPU-based single program, single data (SPSD) parallel programming style to circumvent two major GPU performance bottlenecks, CPU-GPU data transfers and GPU global memory access. In addition, a variable cycling technique is integrated into the method to reduce computational cost when minimizing large-scale nonlinear functions.

## 2. Background and Related Work

This section briefly reviews the background knowledge of interval analysis and existing global optimization methods using interval analysis. The GPU architecture and performance bottlenecks are also summarized in this section. A thorough review of these topics can be found in books written by Moore, Kearfott, and Cloud (26), Hansen and Walster (27), and Hwu, Kirk, and El Hajj (32), among many others.

### *2.1 Global optimization using interval analysis*

Interval analysis, also known as interval arithmetic or interval computation, is a mathematical technique that provides rigorous bounds for the results of numerical computations. In interval analysis, each variable is represented by an interval argument, $[a, b] = \{x: a \le x \le b\}$, where $a$ and $b$ denote the lower and upper bounds for the value of the variable $x$, respectively. A list of interval arithmetic operations is provided in the IEEE 1788 standard (33), and four basic interval arithmetic operations (interval addition, subtraction, multiplication, and division) are defined in Eqs. (3) - (6) for demonstration purposes (26):

$$[x_1, x_2] + [y_1, y_2] = [x_1 + y_1, x_2 + y_2], \tag{3}$$

$$[x_1, x_2] - [y_1, y_2] = [x_1 - y_2, x_2 - y_1], \tag{4}$$

$$[x_1, x_2] \cdot [y_1, y_2] = [\min\{x_1y_1, x_1y_2, x_2y_1, x_2y_2\}, \max\{x_1y_1, x_1y_2, x_2y_1, x_2y_2\}], \tag{5}$$

$$\frac{[x_1, x_2]}{[y_1, y_2]} = [x_1, x_2] \cdot \frac{1}{[y_1, y_2]}. \tag{6}$$

Unlike floating-point arithmetic operations that round the numerical result of each arithmetic operation to its closest binary machine number (25), interval arithmetic operations implement outward rounding to provide rigorous bounds for the result of each arithmetic operation (34). In optimal outward rounding, the outwardly rounded lower bound (left endpoint) is the closest machine number less than or equal to the exact lower bound, and the outwardly rounded upper bound (right endpoint) is the closest machine number greater than or equal to the exact upper bound (26). In practice, the users can also specify their desired significant digits for outward rounding.

Interval evaluation of a mathematical function (i.e., evaluation of a mathematical function through a series of interval arithmetic operations) has been proven to provide rigorous bounds for the value of the mathematical function (26, 27). In other words, when a mathematical function is evaluated with a specified interval argument for each variable, replacing each floating-point arithmetic operation with its corresponding interval arithmetic operation yields a lower bound and an upper bound for the value of the mathematical function, and all possible values of the mathematical function with the specified interval arguments are guaranteed to be included in the range defined by the lower and upper bounds. However, many interval evaluations that include multiple interval arithmetic operations suffer from the dependence problem, also known as interval dependency, and the dependence problem results in overestimated bounds for the value of the mathematical function (i.e., an interval with unwanted excess width). Such overestimated bounds are generated because each occurrence of the variable $x$ is treated as a different variable in the interval evaluation. Several approaches have been introduced to derive exact bounds or at least narrower bounds for the value of a mathematical function when evaluating the function using interval arithmetic operations (26, 27, 35, 36). One major approach is to transform the function into another form, such as a form where every variable only occurs once (if possible), the centered form, the mean value form, or a Taylor series. Another approach, known as splitting or refinement, is to divide the range of each variable in a mathematical function into many subintervals and then take the union of the interval evaluations over these subintervals. It has been proven that the union of these interval evaluations converges to the exact range of the mathematical function when the number of subintervals approaches positive infinity (26).

Interval analysis has been extensively employed in computer-aided proofs in mathematical analysis (26, 37). Several numerical methods using interval analysis (referred to here as existing interval methods) have also been introduced for global optimization (1, 26-30). These existing interval methods often rely on a branch-and-bound process (1, 27, 38, 39), also known as a branch-and-prune process (29, 40) or a branch-and-reduce process (41, 42), to compute the guaranteed bound on the global minimum of the objective function and the corresponding bounds on the value of each variable. Despite several successful applications, these existing interval methods are not commonly adopted for minimizing nonlinear functions in practical applications, especially for minimizing large-scale nonconvex nonlinear functions, because these methods have two major limitations:

- These existing interval methods are not generally applicable to minimizing nonlinear functions. These methods usually employ interval Newton methods and/or consistency techniques (e.g., box consistency and hull consistency (27, 43)) during their iteration process. As a result, these methods usually require the nonlinear function to be continuous and differentiable, and the users often need to transform the nonlinear function into a different form on a case-by-case basis.
- These existing interval methods are designed for CPU-based sequential computing, and they are significantly more computationally expensive than popular numerical optimization methods, such as gradient-based methods and heuristic methods. The recent advances in GPU (31, 32) offer new opportunities to overcome the computational bottleneck in global optimization using interval analysis. Importantly, many personal laptops and desktops are now equipped with a dedicated GPU. Many companies and research institutions are also operating clusters and supercomputers with over 10,000 GPUs (44). However, since GPU is designed for massively parallel computing with a unique compute and memory architecture, these existing interval methods designed for CPU sequential computing cannot utilize the computational power and architecture of GPU.

### *2.2 GPU architecture and performance bottlenecks*

Modern GPUs have brought terascale computing (the ability to perform at least $10^{12}$ calculations per second) to personal laptops and desktops and exascale computing (the ability to perform at least $10^{18}$ calculations per second) to multi-user servers and clusters (31, 32), which enables applications that previously could not be achieved because of long execution times (31). Unlike the CPU, which excels at executing a sequence of instructions, known as a thread, as fast as possible, the GPU excels at executing many threads in parallel and realizes unprecedented computational power through its unique compute and memory architecture (32, 45). The new method introduced in this work is designed and implemented based on GPU compute and memory architecture.

For GPU compute architecture based on the popular Compute Unified Device Architecture (CUDA) (46), threads are organized into a two-level hierarchy, where a grid includes many blocks, and each block includes up to 1024 threads. A modern GPU contains several Streaming Multiprocessors (SMs), and each SM has many cores, also known as processing units or shading units. When threads are assigned to SMs on a GPU, all threads in the same block are assigned to one SM at the same time, and these threads in the same block are further divided into 32-thread units called warps. All threads in a warp always execute the same instruction at any point in time. To efficiently utilize the computational power of a GPU, it is essential to achieve high GPU occupancy by assigning more threads to an SM than the number of cores on the SM (32). In that way, multiple warps can share the same cores on the SM. Specifically, when a warp needs to wait for the result from a long-latency operation, the GPU could select another warp that is ready for execution to use the cores on the SM. In practice, GPU occupancy is calculated as the ratio of the number of active warps on an SM to the maximum number of active warps supported by the SM. To achieve the full occupancy of a GPU, it is necessary to saturate all the SMs by assigning the GPU the theoretical maximum number of active threads that could be supported by the GPU, and the execution of these threads must also satisfy other resource constraints, such as registers and shared memory usages. The theoretical maximum number of active threads could be calculated based on the configuration and compute capability of the GPU. For example, an NVIDIA H100 GPU with the compute capability 9.0 contains 114 SMs, and each SM has 128 cores and supports up to 64 active warps (2,048 threads). To achieve the full occupancy of the NVIDIA H100 GPU, it is necessary to assign at least 233,472 threads (7,296 warps) to the GPU at the same time.

For GPU memory architecture based on CUDA, each thread can access several types of memories on the GPU. As shown in Figure 1, four types of commonly used GPU memories are provided as follows:

- *Global memory*: It provides read-and-write access for all threads in the grid with long access latency and finite access bandwidth. Global memory usually has a storage capacity of several gigabytes (GB).
- *Constant memory*: It provides read-only access for all threads in the grid with short access latency and high access bandwidth. In many modern GPUs, the storage capacity of constant memory is 64 kilobytes (KB).
- *Shared memory*: It provides read-and-write access for threads in the same block and can be accessed at very high speed. Shared memory usually has a storage capacity of tens or hundreds of kilobytes (KB).
- *Registers*: They provide read-and-write access for individual threads and feed directly into the cores on the SM. Each register has a storage capacity of 32 bits.

Modern GPUs have such a unique compute and memory architecture because they are designed for massively parallel computing. As a result, GPU computing is not suitable for all applications. In practice, many applications of GPU computing are based on data parallelism (32, 47), where the computational work to be performed on different portions of the dataset can be done independently. In these applications, each thread uses its corresponding block index and thread index, which are built-in variables in CUDA, to identify the appropriate portion of the dataset to process (32), and all threads execute the same code, known as a kernel function. These applications are instances of the single program, multiple data (SPMD) parallel programming style (32, 48), which belongs to the multiple instruction stream, multiple data stream (MIMD) class based on Flynn's taxonomy (49). Of note, since data parallelism does not exist in existing interval methods for minimizing nonlinear functions, these existing interval methods cannot use the SPMD parallel programming style, and GPU computing has not been commonly adopted for global optimization using interval analysis in practice yet.

Importantly, even in applications that could use the SPMD parallel programming style, GPU computing often suffers from two major performance bottlenecks, CPU-GPU data transfers and GPU global memory access. In typical GPU computing applications using the SPMD parallel programming style, a large amount of data (e.g., several GB) often needs to be transferred between CPU and GPU global memory, and such data transfers are time-consuming. In addition, GPU global memory has long access latency and finite access bandwidth; therefore, it is slow to read and write a large amount of data on GPU global memory. Although several strategies have been introduced to reduce GPU global memory access (32, 50), such as coalescing global memory access and placing reused data on GPU shared memory or registers, it is often difficult to circumvent both of these two major performance bottlenecks, especially CPU-GPU data transfers, in typical GPU computing applications using the SPMD parallel programming style.

## 3. Complete and Rigorous Global Optimization Method Tailored for GPU

To fill the research gap discussed in Section 2, we introduce a GPU-based numerical method for minimizing nonlinear functions. Using interval analysis, coupled with the computational power and architecture of GPU, the method is guaranteed to enclose the global minimum of a nonlinear function, especially a large-scale nonconvex nonlinear function, in a reasonable computation time, even in the presence of rounding errors. Importantly, the method does not endeavor to accelerate any existing interval methods (existing global optimization methods using interval analysis as discussed in Section 2.1) using GPU. Instead, the method is originally designed and implemented for GPU based on GPU compute and memory architecture. This section presents the design and implementation of the GPU-based numerical method. The completeness, rigorousness, and convergence of the method are also discussed in this section.

### *3.1 Method design*

The method is realized by iteratively ruling out the regions in the search domain where the global minimum cannot exist, referred to as the suboptimal regions, and leaving a finite set of regions where the global minimum must exist. A flowchart of the method is illustrated in Figure 2. Each step in the flowchart is described as follows:

*Initialization*: To initialize the iterative process, the *n*-dimensional search domain, also known as the feasible set, is specified by the allowable range of each variable defined in Eq. (2). A list *L* is also created to store the regions (*n*-dimensional "boxes") in the search domain where the global minimum must exist. The list *L* is initialized with a single region that covers the whole search domain.

*Select a region from the list L*: At the beginning of each iteration, among all regions in the list *L*, the region with the smallest lower bound for the value of the nonlinear function in the region is selected and removed from the list *L*, where the lower bound for the value of the nonlinear function in each region in the list *L* is derived by GPU computing in the previous iterations. Of note, that region is selected from the list *L* because it has a better chance of including the global minimum of the nonlinear function than other regions in the list *L* (27).

*Sample over the selected region*: Many sample points are then chosen from the selected region based on a user-specified sampling strategy. Popular sampling techniques, also known as designs for computational experiments (51, 52), such as simple grids, maximum entropy designs, mean squared-error designs, Latin hypercubes, and randomized orthogonal arrays, could be employed to choose sample points in the selected region. Popular gradient-based methods (e.g., interior-point methods and trust region methods) and heuristic methods (e.g., genetic algorithms and simulated annealing) could also be employed to minimize the nonlinear function in the selected region, and their solutions could be used as sample points in the selected region. The nonlinear function is evaluated using interval analysis at each of these sample points, and such interval evaluation yields the lower and upper bounds for the value of the nonlinear function at each sample point. Among all sample points in the selected region, the smallest upper bound for the value of the nonlinear function is used to update the upper bound of the global minimum in the whole search domain, denoted by *GUB*. In other words, the *GUB* value is derived by interval evaluation of the nonlinear function at the best available sample point chosen in all previous iterations and the current iteration, and the global minimum of the nonlinear function is guaranteed to be smaller than or equal to the *GUB* value. In the following steps in the current iteration, the updated *GUB* value is used to rule out the suboptimal regions in the search domain where the global minimum cannot exist. Therefore, if a smaller *GUB* value is efficiently computed through sampling over the selected region, more suboptimal regions could be ruled out in the current iteration. However, it is not necessary to reduce the *GUB* value at a high computational cost (e.g., by choosing a large number of sample points or running a heuristic method with a huge population size in the selected region) in the current iteration because the *GUB* value will be updated in each of the following iterations. In other words, any sampling strategy is acceptable to use in this step, including choosing a single random point in the selected region, but the popular sampling techniques and optimization methods could result in longer computation times in this step, with potential reduction in the total number of iterations, thus reducing the overall computation time of the method.

<u>Remove the suboptimal regions from the list *L*</u>: Using the updated *GUB* value derived from sampling over the selected region, the suboptimal regions in the list *L*, where the global minimum cannot exist, are ruled out and removed from the list *L*. A region in the list *L* is labeled as a suboptimal region and removed from the list *L* if the lower bound for the value of the nonlinear function in the region is larger than the updated *GUB* value. In other words, the value of the nonlinear function in the region is guaranteed to be larger than that at the best sample point found so far, and therefore, the global minimum is impossible to exist in the region.

<u>Partition the selected region and rule out the suboptimal subregions</u>: The selected region is partitioned into many subregions, and the suboptimal subregions in the selected region, where the global minimum cannot exist, are ruled out using the updated *GUB* value. Specifically, the selected region is partitioned into many subregions based on a user-specified partition strategy, such as uniform partition and golden ratio partition. Different scales and coordinate systems could be employed in the partition strategy if necessary. For example, for a 10-dimensional selected region defined by $x_i \in [1, 10^8]$, where $i \in \{1, 2, \ldots, 10\}$, using a uniform partition strategy on a logarithmic scale, the range of each variable in the selected region, $[1, 10^8]$, could be divided into four subintervals, $[1, 10^2]$, $[10^2, 10^4]$, $[10^4, 10^6]$, and $[10^6, 10^8]$, and the selected region is partitioned into $4^{10} = 1{,}048{,}576$ subregions. The nonlinear function is evaluated using interval analysis in all subregions derived from the partition. These interval evaluations yield the lower and upper bounds for the value of the nonlinear function in each subregion. Importantly, as long as the nonlinear function defined in Eq. (1) has a closed-form expression, all possible values of the nonlinear function in a subregion are guaranteed to be included in the range defined by the lower and upper bounds, even if the nonlinear function includes narrow spikes in the subregion. If the lower bound for the value of the nonlinear function in a subregion is larger than the updated *GUB* value, the subregion is labeled as a suboptimal subregion, where the global minimum cannot exist. Importantly, if the selected region is partitioned in all its $n$ dimensions as shown in the example above, both the number of subregions derived from the partition and the GPU computational cost to analyze these subregions increase exponentially with the dimension of the nonlinear function, known as the curse of dimensionality. To circumvent the curse of dimensionality, a variable cycling technique is introduced in Section 3.2, and the integration of the variable cycling technique into the method enables the method to efficiently minimize large-scale nonlinear functions that have more than 100 dimensions.

The two ruling-out processes stated in the last two paragraphs using the updated *GUB* value are generally applicable and do not require the nonlinear function to be continuous or differentiable. When the nonlinear function is a continuous and first-order differentiable function, the first-order derivatives of the nonlinear function can also be used to rule out suboptimal subregions in the selected region. A subregion derived from the partition of the selected region is labeled as a suboptimal subregion if any of the first-order derivatives of the nonlinear function is impossible to be zero in the subregion and the subregion is not located on the edge of the search domain. In other words, a subregion is labeled as a suboptimal region if the necessary first-order optimality conditions for local minimizers (1) are impossible to be satisfied in the subregion. Specifically, a subregion is defined by $x_i \in [\underline{x_i^{(r)}}, \overline{x_i^{(r)}}]$, $i \in \{1, 2, \ldots, n\}$, $r \in \{0, 1, \ldots, k\text{-}1\}$, where $\underline{x_i^{(r)}}$ and $\overline{x_i^{(r)}}$ represent the lower and upper bounds for the value of the variable $x_i$ in the subregion, respectively; $n$ is the dimension of the nonlinear function, and $k$ is the total number of subregions derived from the partition of the selected region. Eq. (2) gives $l_i$ and $u_i$ for each variable and specifies the edges of the search domain. The first-order derivatives of the nonlinear function are evaluated using interval analysis in all subregions derived from the partition of the selected region. These interval evaluations yield the lower and upper bounds for the value of each first-order derivative of the nonlinear function in each subregion, denoted by $DLB_i^{(r)}$ and $DUB_i^{(r)}$. For any $i \in \{1, 2, \ldots, n\}$, if $DLB_i^{(r)} > 0$ and $\underline{x_i^{(r)}} \neq l_i$ in a subregion, the subregion is labeled as a suboptimal subregion. Similarly, for any $i \in \{1, 2, \ldots, n\}$, if $DUB_i^{(r)} < 0$ and $\overline{x_i^{(r)}} \neq u_i$ in a subregion, the subregion is also labeled as a suboptimal subregion.

*Insert the remaining subregions into the list L*: After the suboptimal subregions are ruled out from the selected region based on the updated *GUB* value and the first-order derivatives of the nonlinear function (if available), the remaining subregions (the subregions that are not labeled as suboptimal subregions) in the selected region, where the global minimum may exist, are inserted into the list $L$ for further analysis in the following iterations. Importantly, the two parallel steps shown in the flowchart in Figure 1 only rule out the regions in the search domain where the global minimum is impossible to exist. As a result, the global minimum of the nonlinear function in the search domain is guaranteed to be included in the regions in the list $L$ during the whole iteration process.

*Stopping criteria satisfied*: At the end of each iteration, the method checks whether each user-specified stopping criterion is satisfied. If all stopping criteria are satisfied, the iteration process is halted. Otherwise, a new iteration begins, and among all regions in the updated list $L$, the region with the smallest lower bound for the value of the nonlinear function in the region is selected and removed from the list $L$ for analysis. One or multiple stopping criteria could be specified based on the available computational resources and the practical needs of the user. For example, the stopping criteria could be specified based on the computer run time, the maximum number of iterations, and the region size tolerance (e.g., the interval width in a specific dimension or all dimensions among all the regions in the list $L$).

*Output result*: Once the iteration process is halted, the method outputs all the regions in the list $L$, and the global minimum is guaranteed to be included in these regions. The method also outputs the lower and upper bounds for the value of the global minimum, denoted by *GLB* and *GUB*, and the value of the global minimum is guaranteed to be included in the interval [*GLB*, *GUB*]. The *GLB* value is the smallest lower bound for the value of the nonlinear function in all the regions in the list $L$. The *GUB* value is the updated *GUB* value derived from sampling over the selected region in the last iteration.

### *3.2 Method implementation*

In each iteration of the method, the key step is to partition the selected region into many subregions and rule out suboptimal subregions, where the global minimum cannot exist. This key step is realized through the execution of a kernel function, where each thread in GPU computing evaluates the nonlinear function and the first-order derivatives of the nonlinear function (if available) in one corresponding subregion using interval analysis independently. To efficiently utilize the computational power of the user's GPU(s), the kernel function is implemented based on GPU compute and memory architecture. In addition, a variable cycling technique is integrated into the kernel function for minimizing large-scale nonlinear functions that have more than 100 dimensions. Moreover, a data representation with limited memory usage is employed for each region in the list $L$ to avoid host RAM overflow.

To achieve high GPU occupancy based on GPU compute architecture, the number of subregions derived from the partition of the selected region is determined based on the configuration and compute capability of the user's GPU(s). Specifically, since each subregion corresponds to one thread in GPU computing, it is necessary to derive an appropriate number of subregions from the partition to saturate the SMs of the user's GPU(s). In practice, the theoretical maximum number of active threads that could be supported by the user's GPU(s) serves as a reference value for the number of subregions. For example, as discussed in Section 2.2, for a server with eight NVIDIA H100 GPUs, at least 1,867,776 subregions need to be derived from the partition of the selected region to achieve the full occupancy of these eight GPUs. If the number of subregions derived from the partition is significantly larger (e.g., more than two orders of magnitude larger) than the reference value, a large portion of the threads will be executed in sequential computing rather than in parallel computing, leading to a long execution time. If the number of subregions derived from the partition is significantly smaller (e.g., less than the total number of cores on the user's GPU(s)) than the reference value, many GPU cores will become idle during the execution of the kernel function, leading to GPU underutilization. In addition, since 32 threads in a warp are always scheduled to execute the same instruction at the same time in GPU computing, it is beneficial to choose the number of subregions derived from the partition as a multiple of 32.

To circumvent two major GPU performance bottlenecks (CPU-GPU data transfers and GPU global memory access) based on GPU memory architecture, a novel GPU-based single program, single data (SPSD) parallel programming style is employed in the implementation of the kernel function. Specifically, to evaluate the nonlinear function and the first-order derivatives of the nonlinear function (if available) in each subregion using interval analysis, the data set that each thread in GPU computing needs to use is the location of its corresponding subregion in the search domain. The location of a subregion in the search domain is represented by the lower and upper bounds for the value of each variable in the subregion, denoted by $\underline{x_i^{(r)}}$ and $\overline{x_i^{(r)}}$, respectively, where $x_i \in [\underline{x_i^{(r)}}, \overline{x_i^{(r)}}]$, $i \in \{1, 2, \ldots, n\}$, $r \in \{0, 1, \ldots, k\text{-}1\}$ in the subregion. If the common single program, multiple data (SPMD) parallel programming style discussed in Section 2.2 is employed, the selected region needs to be partitioned into many subregions on CPU, and the location data of these subregions in the search domain ($2 \times n \times k$ floating-point numbers) are transferred to GPU global memory. Each thread in GPU computing identifies its corresponding subregion using its block index and thread index, and each thread reads its corresponding subregion location data from GPU global memory. Importantly, the implementation of the kernel function based on the common SPMD parallel programming style suffers from the two major GPU performance bottlenecks since it requires a large amount of data transfer from CPU to GPU and many GPU global memory reads at the beginning of the kernel function execution. To circumvent these GPU performance bottlenecks, a GPU-based single program, single data (SPSD) parallel programming style is introduced and employed. Using the GPU-based SPSD parallel programming style, only the location of the selected region ($2 \times n$ floating-point numbers) is transferred to GPU constant memory. The selected region is partitioned on GPU, and each thread *computes* its corresponding subregion location using its block index, thread index, and the location of the selected region. As an example, when minimizing a two-dimensional nonlinear function ($n = 2$), the selected region in the current iteration is defined by $x_1 \in [\underline{x_1^{select}}, \overline{x_1^{select}}]$ and $x_2 \in [\underline{x_2^{select}}, \overline{x_2^{select}}]$, and the range of each variable in the selected region is uniformly divided into four subintervals, resulting in 16 subregions in total ($k = 16$). As illustrated in Figure 3(a), each subregion is assigned a subregion index, denoted by Sidx, and Sidx $\in \{0, 1, \ldots, 15\}$ in this example. Each thread in GPU computing has its unique index in the grid, also known as the absolute position of the thread in the entire grid of blocks, and the unique index $r$ of each thread is derived by

$$r = \text{threadIdx.x} + \text{blockIdx.x} \times \text{blockDim.x}, \qquad (7)$$

where threadIdx.x is the thread index, blockIdx.x is the block index, blockDim.x is the block size (number of threads in a block), and one-dimensional grid and blocks are used in this example. The thread index and block index are two built-in variables in CUDA, and the block size is a constant specified by the user when the kernel function is launched. The unique index of each thread in the grid is mapped to the subregion index. In other words, a thread with the unique index of $r$ in the grid corresponds to the subregion with Sidx = $r$. As shown in Figure 3(b), the subinterval indices of each subregion, denoted by $I_1^{(r)}$ and $I_2^{(r)}$ in this example, are computed by

$$I_1^{(r)} = r \,\%\, 4, \tag{8}$$

$$I_2^{(r)} = r \,//\, 4, \tag{9}$$

where % in Eq. (8) represents a modulo operation that returns the remainder of a division, and // in Eq. (9) represents a floor division operation that performs a division and rounds the result down to the nearest integer. The location of a subregion in the search domain, represented by $\underline{x_i^{(r)}}$ and $\overline{x_i^{(r)}}$, is computed using the location of the selected region, represented by $\underline{x_i^{select}}$ and $\overline{x_i^{select}}$, and the subinterval indices of the subregion:

$$\underline{x_i^{(r)}} = \underline{x_i^{select}} + \frac{\overline{x_i^{select}} - \underline{x_i^{select}}}{4} \cdot I_i^{(r)}, \tag{10}$$

$$\overline{x_i^{(r)}} = \underline{x_i^{select}} + \frac{\overline{x_i^{select}} - \underline{x_i^{select}}}{4} \cdot \left( I_i^{(r)} + 1 \right), \tag{11}$$

where $i \in \{1, 2\}$ in this example. Of note, the example above with $n = 2$ and $k = 16$ is provided here for demonstration purposes only. In practice, the selected region is usually partitioned into at least thousands of subregions ($k > 1{,}000$) to achieve high GPU occupancy, and the subinterval indices of each subregion are computed through a series of modulo and floor division operations when the nonlinear function has more than two dimensions ($n > 2$). Importantly, using the GPU-based SPSD parallel programming style, only limited data ($2 \times n$ floating-point numbers) needs to be transferred from CPU to GPU constant memory, and all threads in GPU computing read the same data from GPU constant memory. When minimizing large-scale nonlinear functions where GPU constant memory does not have enough space to store the location data of the selected region (e.g., $n > 4{,}000$ with FP64 format), such location data could be transferred from CPU to GPU global memory. Since all threads in GPU computing read the same location data saved on GPU global memory, the data will be automatically cached on GPU, leading to minimal GPU global memory reads. After the execution of the kernel function, as stated in Section 3.1, the subregions that are not labeled as suboptimal subregions (referred to as the remaining subregions) are transferred from GPU to CPU and inserted into the list $L$. If there are many remaining subregions and each remaining subregion is represented by its location data, a large amount of data needs to be written onto GPU global memory and transferred to CPU ($2 \times n$ floating-point numbers for each remaining subregion). To reduce GPU global memory writes and the amount of data that needs to be transferred from GPU to CPU, the index of each remaining subregion, Sidx, is written in an array on GPU global memory (one integer for each remaining subregion), and the array is transferred to CPU.

When minimizing large-scale nonlinear functions, the selected region in the search domain is partitioned using a variable cycling technique, where the selected region is partitioned in a portion of its $n$ dimensions in each iteration in a cyclic manner, significantly reducing the computational cost during the iteration process. As discussed in Section 3.1, if the selected region is partitioned in all its $n$ dimensions in the kernel function, the number of subregions and the corresponding GPU computational cost increase exponentially with the dimension of the nonlinear function. To reduce the number of subregions that need to be analyzed during the execution of the kernel function, the selected region is partitioned in a portion of its $n$ dimensions based on the cycling index of the selected region. For example, when minimizing a nonlinear function with $n = 1{,}000$, 10 out of 1,000 dimensions of the selected region are partitioned in each iteration, and the range of each of these 10 variables is equally divided into four subintervals, while the ranges of the other $n$-10 variables are intact, resulting in $4^{10} = 1{,}048{,}576$ subregions from the partition of the selected region in each iteration, aligned with the configuration and compute capability of a modern GPU for demonstration purposes. At the first iteration, the single region that covers the whole search domain is assigned to a cycling index of 1, and the first 10 dimensions ($x_1$ to $x_{10}$) of the selected region are partitioned. The remaining subregions derived from the first iteration, where the global minimum must exist, are assigned to a cycling index of 11 when they are inserted into the list $L$, and these regions will be partitioned in the next 10 dimensions ($x_{11}$ to $x_{20}$) if any of them is selected in the following iterations, and so on. If the selected region in an iteration has a cycling index of 991, the last 10 dimensions ($x_{991}$ to $x_{1000}$) of the selected region are partitioned in the iteration. The remaining subregions derived from the iteration are assigned to a cycling index of 1 when they are inserted into the list $L$, and these regions will be partitioned in the first 10 dimensions ($x_1$ to $x_{10}$) if any of them is selected in the following iterations. Importantly, using the variable cycling technique, a smaller number of regions needs to be partitioned during the whole iteration process, because many suboptimal regions where the global minimum cannot exist are ruled out before they are partitioned in all their $n$ dimensions, significantly reducing the computational cost of the method when minimizing large-scale nonlinear functions.

During the iteration process of the method, the list $L$ may contain a large number of regions (e.g., more than $1\times10^8$ regions) where the global minimum must exist. To avoid host RAM overflow, only the following four types of data are stored in the host RAM for each region in the list $L$:

- *Subregion index*: It is the index of the subregion that is derived from the partition of the selected region, as shown in Figure 3(a), denoted by Sidx, and Sidx $\in \{0, 1, \ldots, k\text{-}1\}$, where $k$ is the total number of subregions derived from the partition of the selected region. The subregion index is stored in the host RAM as one integer.
- *Iteration index*: It is the index of the iteration in which the subregion is derived from the partition of the selected region and inserted into the list $L$, denoted by itr, and itr $\in \{0, 1, \ldots, t\text{-}1\}$, where $t$ is the total number of iterations of the method. The iteration index is stored in the host RAM as one integer.
- *Cycling index*: It is the index that determines the specific dimensions of the region to be partitioned if the region is selected at the beginning of an iteration. The cycling index is stored in the host RAM as one integer.
- *Lower bound for the value of the nonlinear function in the region*: It is computed by the interval evaluation of the nonlinear function in the subregion using GPU computing. The lower bound for the value of the nonlinear function in the region is stored in the host RAM as one floating-point number.

Other types of data, such as the upper bound for the value of the nonlinear function in each region and the location data of each region in the list $L$, do not need to be stored in the host RAM. The list $L$ with these four types of data could be realized through a dictionary-like data structure or multiple arrays. When a region is selected from the list $L$ at the beginning of an iteration, the location data of the selected region, represented by the lower and upper bounds for the value of each variable in the selected region ($2\times n$ floating-point numbers), is computed through the process exemplified by Figure 3 and Eqs. (8) - (11) on CPU, and the location data of the selected region in each iteration is stored in the host RAM.

### *3.3 Method completeness, rigorousness, and convergence*

Global optimization methods are classified based on their completeness and rigorousness (1):

- A *complete* method reaches the global minimum with certainty and encloses the global minimum within prescribed tolerances after a finite computation time.
- A *rigorous* method encloses the global minimum within prescribed tolerances even in the presence of rounding errors, except in near-degenerate cases.

According to the definitions above, the method introduced in this section is a complete and rigorous global optimization method. For the method completeness, since the method only rules out the regions in the search domain where the global minimum is guaranteed not to exist, the method never misses the global minimum, and the global minimum must exist in the regions outputted by the method. For the method rigorousness, the method employs interval analysis to take rounding errors into account and provides safe bounds for both the location and the value of the global minimum.

When the global minimum is not located at positive or negative infinity and a user-specified region size tolerance is used as the stopping criterion, the method is guaranteed to converge after a finite computation time because the regions in the list $L$ are iteratively selected and partitioned, and the sizes of these regions keep decreasing during the iteration process of the method. The method converges quickly when the nonlinear function is separable or does not suffer from a severe dependence problem, as demonstrated by minimizing the 10 multimodal benchmark test functions in Section 4, where the computation time of the method increases approximately quadratically with the dimension of these 10 multimodal benchmark test functions. When the nonlinear function suffers from a severe dependence problem, or the global minimum of the nonlinear function resides in a flat valley, the method converges more slowly but the computation time does not increase exponentially with the dimension of the nonlinear function, as demonstrated by minimizing the Rosenbrock function in Section 4, where the computational time of the method has a cubic to quartic increase with the dimension of the Rosenbrock function.

## 4. Computational Experiments

The GPU-based numerical method introduced in Section 3 is validated by minimizing 11 benchmark test functions with scalable dimensions, including the well-known Ackley function, Griewank function, Levy function, Rastrigin function, and Rosenbrock function. The method successfully encloses the guaranteed global minimum of these 11 benchmark test functions with up to 10,000 dimensions using only one GPU in a reasonable computation time. In addition, when minimizing the Ackley function with 100 dimensions, the results of the method are compared with those of seven popular numerical optimization methods. All seven popular numerical optimization methods fail to find the global minimum of the Ackley function, while the method successfully encloses the guaranteed global minimum of the Ackley function. Moreover, the method successfully encloses the guaranteed global minimum of a discontinuous test function with scalable dimensions built from the Levy function and the Dirac delta function, demonstrating the ability of the method to minimize discontinuous nonlinear functions.

The 10 multimodal benchmark test functions with scalable dimensions are defined in SI Appendix A. Each of these 10 nonlinear functions is highly nonconvex and has many local minima. These benchmark test functions represent grand challenges of global optimization, and enclosing the guaranteed global minimum of these benchmark test functions with more than 80 dimensions has not been reported in the literature, where enclosing the *guaranteed* global minimum requires the optimization method to be complete and rigorous as defined in Section 3.3.

To validate the GPU-based numerical method introduced in Section 3 and evaluate the computational costs of the method when minimizing nonlinear functions with different dimensions, the method is employed to minimize each of these 10 benchmark test functions with 50, 100, 500, and 1,000 dimensions ($n$ = 50, 100, 500, and 1,000 in these functions) using a laptop with one GPU, a workstation with one GPU, a local server with one GPU, and a cloud server with one GPU, and the configurations of these four typical computational resources are provided in SI Appendix B. The method is also employed to minimize the Levy function with 2,000, 5,000, and 10,000 dimensions ($n$ = 2,000, 5,000, and 10,000 in the Levy function) using the cloud server with one GPU to demonstrate its ability to minimize large-scale nonlinear functions.

When minimizing these 10 multimodal benchmark test functions with different dimensions using the method, double-precision floating-point format (FP64 format) is employed for CPU and GPU computing during the iteration process of the method, and the same sampling strategy, partition strategy, and stopping criterion are employed in the method for demonstration purposes. For the step to sample over the selected region in each iteration, 10 out of $n$ dimensions of the selected region are partitioned using the variable cycling technique introduced in Section 3.2, and the range of each of these 10 variables is equally divided into four subintervals, while the ranges of the other $n$-10 variables are intact. The partition of the selected region results in $4^{10}$ = 1,048,576 subregions, and each subregion is denoted by $x_i \in [\underline{x_i^{(r)}}, \overline{x_i^{(r)}}]$, where $i \in \{1, 2, \ldots, n\}$, $r \in \{0, 1, \ldots, 1{,}048{,}575\}$. 10 sample points are chosen along the diagonal of each subregion, where the coordinates of the two diagonal endpoints are $[\underline{x_1^{(r)}}, \underline{x_2^{(r)}}, \ldots, \underline{x_n^{(r)}}]$ and $[\overline{x_1^{(r)}}, \overline{x_2^{(r)}}, \ldots, \overline{x_n^{(r)}}]$, and the 10 sample points equally divide the diagonal into 11 pieces. As detailed in Section 3.1, the *GUB* value is updated based on the interval evaluations of the benchmark test function at these sample points using GPU computing. For the step to partition the selected region and rule out the suboptimal subregions in each iteration, the same partition strategy as the sampling step stated in this paragraph is employed. Among the 1,048,576 subregions derived from the partition of the selected region, the suboptimal subregions are ruled out based on the updated *GUB* value and the first-order derivatives of the benchmark test function, as detailed in Section 3.1. The region width tolerance of $1\times10^{-4}$ is employed as the stopping criterion. Specifically, at the end of each iteration, the iteration process is halted if the widths in all $n$ dimensions of all regions in the list $L$ are less than $1\times10^{-4}$.

The results of these computational experiments show that the GPU-based numerical method successfully encloses the known global minimum for each of these 10 benchmark test functions with different dimensions in a reasonable computation time, demonstrating the effectiveness of the method. The results of these 10 benchmark test functions with 1,000 dimensions using the workstation, local server, and cloud server are summarized in Table 1. As shown in Table 1, in each computational experiment, the method outputs only one region that satisfies the region width tolerance of $1\times10^{-4}$ and includes the known global minimum of the benchmark test function. In addition, using different computational resources to minimize a benchmark test function with a specific dimension, the method with the same sampling strategy, partition strategy, and stopping criterion as stated in the last paragraph always outputs the same region, because random variation is not involved in the method.

The results of these computational experiments also show that the computational time of the method increases approximately quadratically with the dimension of the benchmark test function, demonstrating the efficiency of the method, especially the variable cycling technique introduced in Section 3.2, to minimize large-scale nonlinear functions. As an example, the results of the Levy function with different dimensions using the laptop, workstation, local server, and cloud server are summarized in Table 2. The results of the other nine benchmark test functions with different dimensions are provided in SI Appendix C. As shown in these tables, the number of iterations increases linearly with the dimension of the benchmark test function. Importantly, the number of terms in each benchmark test function increases linearly with the dimension of the function, and a major part of the computation time in each iteration of the method is spent on the interval evaluation of the benchmark test function. As a result, the computational time of each iteration increases approximately linearly with the dimension of the benchmark test function, leading to an approximately quadratic increase in the total computation time with the dimension of the benchmark test function.

Besides the 10 multimodal benchmark test functions, the method is employed to minimize the Rosenbrock function with scalable dimensions using the laptop, workstation, local server, and cloud server (configurations provided in SI Appendix B). The Rosenbrock function with scalable dimensions is defined in SI Appendix D. Although the Rosenbrock function does not have many local minima, the function is highly nonconvex and has strong variable coupling, leading to a severe dependence problem. Moreover, the global minimum of the Rosenbrock function resides in a long, narrow, parabolic-shaped flat valley, making it difficult for optimization methods to converge to the global minimum. When minimizing the Rosenbrock function using the method, the same sampling strategy, partition strategy, and stopping criterion as in minimizing the 10 multimodal benchmark test functions are employed for demonstration purposes. The results of this computational experiment show that the GPU-based numerical method successfully encloses the known global minimum for the Rosenbrock function with 50, 100, 500, and 1,000 dimensions ($n$ = 50, 100, 500, and 1,000 in the function) using the four different computational resources. As shown in the table provided in SI Appendix D, the method outputs several interconnected regions that include the known global minimum of the function, and the computational time of the method has a cubic to quartic (but not exponential) increase with the dimension of the function.

To compare the GPU-based numerical method with popular CPU-based numerical optimization methods, especially popular global optimization methods based on regular floating-point arithmetic, seven popular CPU-based numerical optimization methods (basin-hopping, BFGS, differential evolution, DIRECT, interior-point, Nelder-Mead, and trust region), implemented in a commercial software, MATLAB, and a Python package, SciPy, are employed to minimize the Ackley function with 100 dimensions as defined in SI Appendix A using the laptop and the local server (configurations provided in SI Appendix B). Among these seven popular numerical optimization methods, the DIRECT method does not require any starting point to initialize the method and always produces the same solution in multiple runs, while the other six methods require starting points to initialize these methods and often produce different solutions with various starting points. As a result, each of these six methods is run multiple times with random starting points in the search domain. The results of this computational experiment using the laptop and the local server are provided in SI Appendix E. As shown in these results, all seven popular numerical optimization methods fail to find the known global minimum of the Ackley function, even with multiple runs and a longer total computation time than the GPU-based numerical method. In contrast, the GPU-based numerical method successfully encloses the known global minimum in a single run.

To demonstrate the ability of the GPU-based numerical method to minimize discontinuous nonlinear functions, the method is employed to minimize a discontinuous test function with scalable dimensions using the laptop, workstation, local server, and cloud server (configurations provided in SI Appendix B). As defined in SI Appendix F, the discontinuous test function is built from the Levy function and the Dirac delta function, and the global minimum is introduced by the Dirac delta function. When minimizing the discontinuous test function using the method, the same sampling strategy, partition strategy, and stopping criterion as in minimizing the 10 multimodal benchmark test functions are employed for demonstration purposes. The results of this computational experiment show that the GPU-based numerical method successfully encloses the known global minimum for the discontinuous test function with 50, 100, 500, and 1,000 dimensions ($n$ = 50, 100, 500, and 1,000 in the function) using the four different computational resources. As shown in the table provided in SI Appendix F, the method outputs only one region that includes the known global minimum of the discontinuous test function, and the computational time of the method increases approximately quadratically with the dimension of the discontinuous test function.

## 5. Final Remarks

We introduce a GPU-based numerical method to minimize nonlinear functions subject to simple bounds on the variables. The method is realized through an iterative process that rules out the regions in the search domain where the global minimum cannot exist and leaves a finite set of regions where the global minimum must exist. The iteration process is designed and implemented based on GPU compute and memory architecture, where a GPU-based single program, single data parallel programming style is employed to circumvent major GPU performance bottlenecks, and a variable cycling technique is employed to reduce computational cost when minimizing large-scale nonlinear functions. The computational experiments on 11 multimodal benchmark test functions with scalable dimensions show that our method successfully encloses the global minimum of these benchmark test functions with up to 10,000 dimensions using only one GPU in a reasonable computation time.

We project that the GPU-based numerical method can enclose the guaranteed global minimum of large-scale nonconvex nonlinear functions that have more than one million dimensions and many local minima using a GPU cluster in many practical applications. In addition, various sampling and partition strategies could be explored to further improve the efficiency of the method. Moreover, the method can be extended to solve optimization problems subject to equality and inequality constraints.

**Data, Materials, and Software Availability**

All data are included in the manuscript and SI Appendices. The source code of computational experiments is available on GitHub (https://github.com/CMU-Integrated-Design-Innovation-Group/GPUGO).

**Acknowledgments**

The authors gratefully acknowledge the support of Lambda, Inc. for providing the cloud server with NVIDIA GH200 Grace Hopper Superchip for computational experiments in this work.

**Funding Statement**

The authors received no specific funding for this work.

**Figures and Tables**

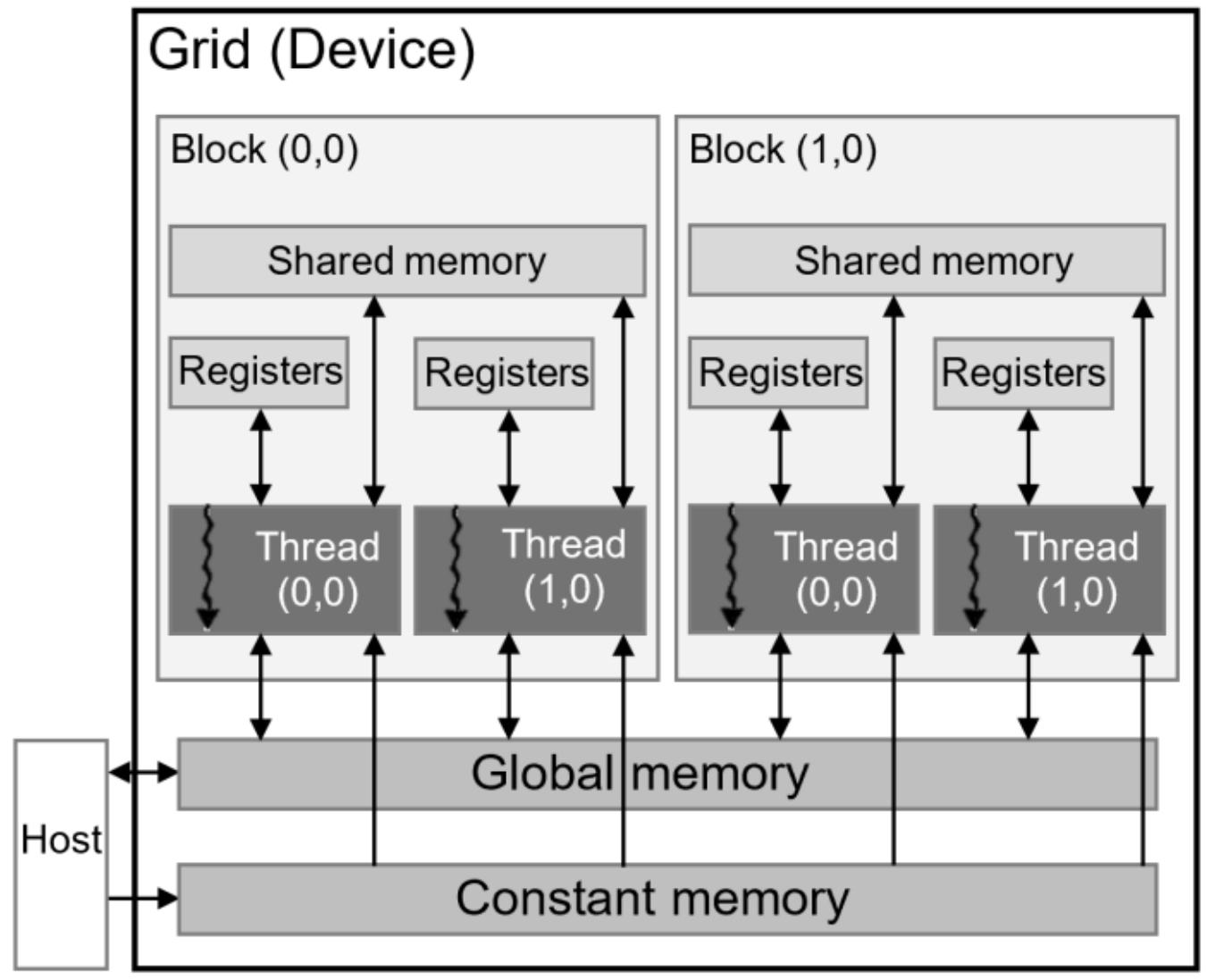

**Figure 1.** An overview of GPU memory architecture based on CUDA (32)

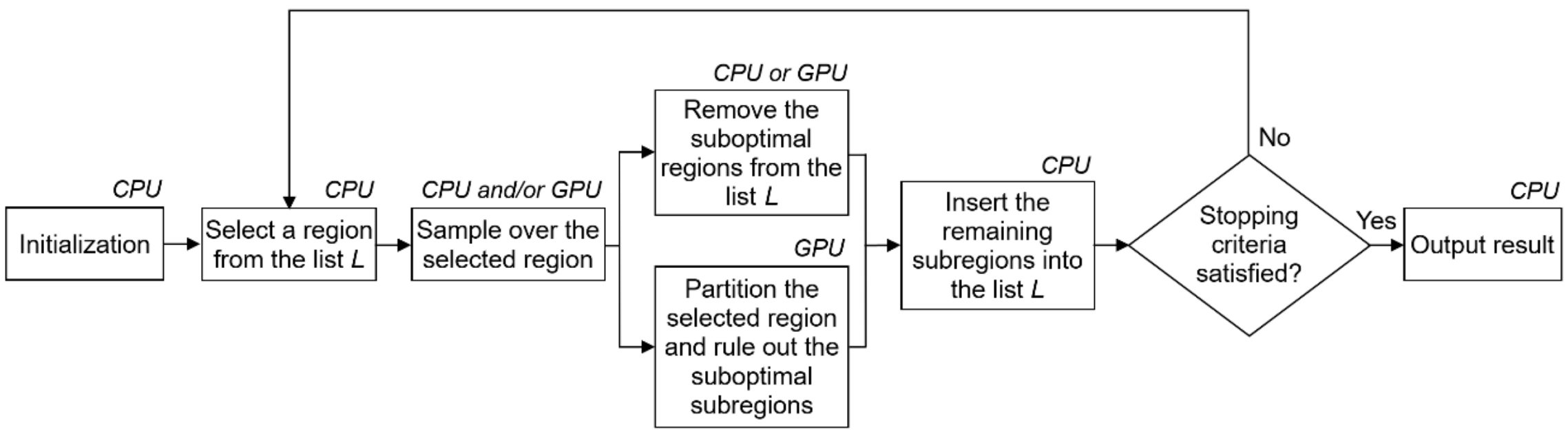

**Figure 2.** The flowchart of the GPU-based numerical method

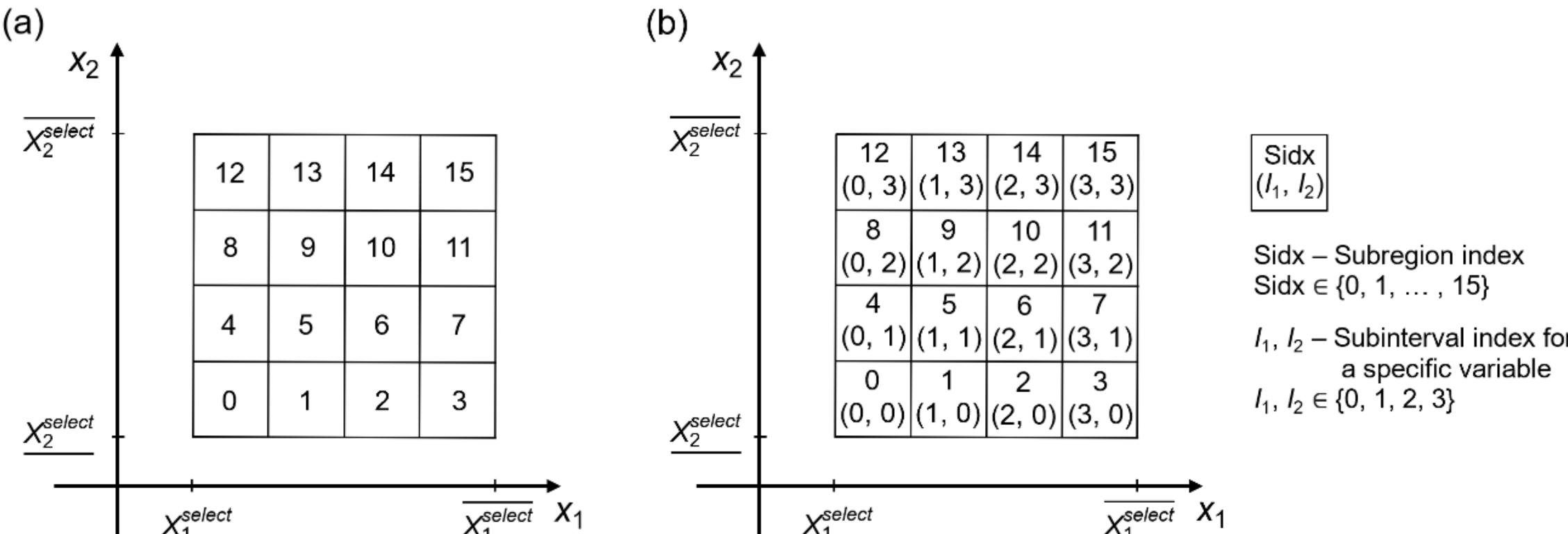

**Figure 3.** Partition of the selected region and indexing of each subregion

**Table 1.** Results of 10 multimodal benchmark test functions with 1,000 dimensions (*n* = 1,000)

| Function | Number of Iterations | Number of Output Regions | Computation Time | | |
|---|---|---|---|---|---|
| | | | Workstation | Local Server | Cloud Server |
| Ackley | 1,000 | 1 | 33,025 s (9.17 h) | 72,625 s (20.17 h) | 2,288 s (0.64 h) |
| Belegundu | 900 | 1 | 23,780 s (6.61 h) | 37,683 s (10.47 h) | 1,896 s (0.53 h) |
| Breiman | 800 | 1 | 17,041 s (4.73 h) | 53,621 s (14.89 h) | 1,218 s (0.34 h) |
| Fu | 900 | 1 | 55,774 s (15.49 h) | 148,104 s (41.14 h) | 3,538 s (0.98 h) |
| Griewank | 1,100 | 1 | 37,060 s (10.29 h) | 82,044 s (22.79 h) | 2,523 s (0.70 h) |
| Levy | 900 | 1 | 52,927 s (14.70 h) | 117,176 s (32.55 h) | 3,582 s (1.00 h) |
| Rastrigin | 900 | 1 | 23,214 s (6.45 h) | 71,765 s (19.93 h) | 1,727 s (0.48 h) |
| Salomon | 1,100 | 1 | 13,606 s (3.78 h) | 31,409 s (8.72 h) | 1,435 s (0.40 h) |
| Styblinski | 900 | 1 | 24,218 s (6.73 h) | 64,924 s (18.03 h) | 1,791 s (0.50 h) |
| Zabinsky | 800 | 1 | 36,699 s (10.19 h) | 80,920 s (22.48 h) | 2,635 s (0.73 h) |

Our method outputs only one region that satisfies the region width tolerance of $1\times10^{-4}$ and includes the known global minimum for each of these 10 benchmark test functions with 1,000 dimensions.

**Table 2.** Results of the Levy function with different dimensions

| Dimension | Number of Iterations | Number of Output Regions | Computation Time | | | |
|---|---|---|---|---|---|---|
| | | | Laptop | Workstation | Local Server | Cloud Server |
| *n* = 50 | 45 | 1 | 644 s | 136 s | 311 s | 15 s |
| *n* = 100 | 90 | 1 | 2,537 s (0.70 h) | 536 s (0.15 h) | 1,184 s (0.33 h) | 42 s (0.01 h) |
| *n* = 500 | 450 | 1 | 63,063 s (17.52 h) | 13,258 s (3.68 h) | 29,217 s (8.12 h) | 909 s (0.25 h) |
| *n* = 1,000 | 900 | 1 | NA | 52,927 s (14.70 h) | 117,176 s (32.55 h) | 3,582 s (1.00 h) |
| *n* = 2,000 | 1,800 | 1 | NA | 215,458 s (59.85 h) | NA | 14,543 s (4.04 h) |
| *n* = 5,000 | 4,500 | 1 | NA | NA | NA | 93,447 s (25.96 h) |
| *n* = 10,000 | 9,000 | 1 | NA | NA | NA | 376,589 s (104.61 h) |

Our method outputs only one region that satisfies the region width tolerance of $1 \times 10^{-4}$ and includes the known global minimum for the Levy function with different dimensions. The number of iterations increases linearly with the dimension of the Levy function. The computational time of our method increases approximately quadratically with the dimension of the Levy function. “NA” in the table represents that the corresponding computational resource is not used to minimize the Levy function with the specific dimension because of an expected long execution time.

**Supporting Information for**

# Global optimization tailored for graphics processing units: Complete and rigorous search for large-scale nonlinear minimization

Guanglu Zhang[a*], Qihang Shan[a], Jonathan Cagan[a]
[a]Department of Mechanical Engineering, Carnegie Mellon University, Pittsburgh, PA 15213

*Corresponding author: Guanglu Zhang
Email: glzhang@cmu.edu

**This PDF file includes:**

Supporting text
Figure S1
Tables S1 to S14
SI References

## Appendix A. The 10 multimodal benchmark test functions with scalable dimensions

The 10 multimodal benchmark test functions with scalable dimensions are defined below.

The Ackley function has the form (1-3):

$$f(\boldsymbol{x}) = -20\exp\left(-0.02\sqrt{\frac{1}{n}\sum_{i=1}^{n} x_i^2}\right) - \exp\left(\frac{1}{n}\sum_{i=1}^{n}\cos(2\pi x_i)\right) + 20 + \mathrm{e}, \tag{S1}$$

$$-35 \le x_i \le 40, \tag{S2}$$

where $i \in \{1, 2, \ldots, n\}$. The global minimum of the Ackley function is located at $\boldsymbol{x}^* = [0, 0, \ldots, 0]$, and $f(\boldsymbol{x}^*) = 0$.

The Belegundu function, also known as the deflected corrugated spring function, has the form (4):

$$f(\boldsymbol{x}) = 0.1\sum_{i=1}^{n}(x_i - 5)^2 - \cos\left[5\sqrt{\sum_{i=1}^{n}(x_i - 5)^2}\right], \tag{S3}$$

$$-10 \le x_i \le 11, \tag{S4}$$

where $i \in \{1, 2, \ldots, n\}$. The global minimum of the Belegundu function is located at $\boldsymbol{x}^* = [5, 5, \ldots, 5]$, and $f(\boldsymbol{x}^*) = -1$.

The Breiman function, also known as the cosine mixture function, has the form (5):

$$f(\boldsymbol{x}) = -0.1\sum_{i=1}^{n}\cos(5\pi x_i) + \sum_{i=1}^{n} x_i^2, \tag{S5}$$

$$-1 \le x_i \le 2, \tag{S6}$$

where $i \in \{1, 2, \ldots, n\}$. The global minimum of the Breiman function is located at $\boldsymbol{x}^* = [0, 0, \ldots, 0]$, and $f(\boldsymbol{x}^*) = -0.1n$.

The Fu function, also known as the trigonometric function, has the form (6):

$$f(\boldsymbol{x}) = 1 + \sum_{i=1}^{n} 8\left\{\sin\left[7(x_i - 0.9)^2\right]\right\}^2 + 6\left\{\sin\left[14(x_i - 0.9)^2\right]\right\}^2 + (x_i - 0.9)^2, \tag{S7}$$

$$-10 \le x_i \le 10, \tag{S8}$$

where $i \in \{1, 2, \ldots, n\}$. The global minimum of the Fu function is located at $\boldsymbol{x}^* = [0.9, 0.9, \ldots, 0.9]$, and $f(\boldsymbol{x}^*) = 1$.

The Griewank function has the form (3, 7):

$$f(\boldsymbol{x}) = 1 + \sum_{i=1}^{n}\frac{x_i^2}{4000} - \prod_{i=1}^{n}\cos\left(\frac{x_i}{\sqrt{i}}\right), \tag{S9}$$

$$-100 \le x_i \le 110, \tag{S10}$$

where $i \in \{1, 2, \ldots, n\}$. The global minimum of the Griewank function is located at $\boldsymbol{x}^* = [0, 0, \ldots, 0]$, and $f(\boldsymbol{x}^*) = 0$.

The Levy function has the form (8):

$$f(\boldsymbol{x}) = \frac{\pi}{n}\left\{10\sin^2(\pi y_1) + (y_n - 1)^2 + \sum_{i=1}^{n-1}\left[(y_i - 1)^2\left(1 + 10\sin^2(\pi y_{i+1})\right)\right]\right\}, \quad \text{(S11)}$$

$$y_i = 1 + 0.25(x_i - 1), \quad \text{(S12)}$$

$$-10 \leq x_i \leq 10, \quad \text{(S13)}$$

where $i \in \{1, 2, \ldots, n\}$. The global minimum of the Levy function is located at $\boldsymbol{x}^* = [1, 1, \ldots, 1]$, and $f(\boldsymbol{x}^*) = 0$.

The Rastrigin function has the form (9, 10):

$$f(\boldsymbol{x}) = 10n + \sum_{i=1}^{n}\left[x_i^2 - 10\cos(2\pi x_i)\right], \quad \text{(S14)}$$

$$-5.5 \leq x_i \leq 6, \quad \text{(S15)}$$

where $i \in \{1, 2, \ldots, n\}$. The global minimum of the Rastrigin function is located at $\boldsymbol{x}^* = [0, 0, \ldots, 0]$, and $f(\boldsymbol{x}^*) = 0$.

The Salomon function has the form (11):

$$f(\boldsymbol{x}) = 1 - \cos\left(2\pi\sqrt{\sum_{i=1}^{n} x_i^2}\right) + 0.1\sqrt{\sum_{i=1}^{n} x_i^2}, \quad \text{(S16)}$$

$$-100 \leq x_i \leq 110, \quad \text{(S17)}$$

where $i \in \{1, 2, \ldots, n\}$. The global minimum of the Salomon function is located at $\boldsymbol{x}^* = [0, 0, \ldots, 0]$, and $f(\boldsymbol{x}^*) = 0$.

The Styblinski function, also known as the multimodal corrupted convex quadratic function, has the form (12):

$$f(\boldsymbol{x}) = \frac{1}{2n}\sum_{i=1}^{n} x_i^2 - 4n\prod_{i=1}^{n}\cos(x_i), \quad \text{(S18)}$$

$$-10 \leq x_i \leq 11, \quad \text{(S19)}$$

where $i \in \{1, 2, \ldots, n\}$. The global minimum of the Styblinski function is located at $\boldsymbol{x}^* = [0, 0, \ldots, 0]$, and $f(\boldsymbol{x}^*) = -4n$.

The Zabinsky function, also known as the sinusoidal function, has the form (13, 14):

$$f(\boldsymbol{x}) = -2.5\prod_{i=1}^{n}\sin\left(x_i - \frac{\pi}{6}\right) - \prod_{i=1}^{n}\sin\left[5(x_i - \frac{\pi}{6})\right], \quad \text{(S20)}$$

$$0 \leq x_i \leq \pi, \quad \text{(S21)}$$

where $i \in \{1, 2, \ldots, n\}$. The global minimum of the Zabinsky function is located at $\boldsymbol{x}^* = [2\pi/3, 2\pi/3, \ldots, 2\pi/3]$, and $f(\boldsymbol{x}^*) = -3.5$.

## Appendix B. Configurations of four typical computational resources

The configurations of four typical computational resources, including a laptop, a workstation, a local server, and a cloud server, are provided in Table S1. Only one GPU in each computational resource is used in the computational experiments detailed in Section 4.

**Table S1.** Configurations of four typical computational resources

| Resources | CPU | RAM | GPU |
|---|---|---|---|
| Laptop | Intel Core 12th Gen i7-12700H | 16 GB | NVIDIA GeForce RTX 3070Ti |
| Workstation | AMD Ryzen Threadripper PRO 7975WX | 256 GB | NVIDIA GeForce RTX 4090 |
| Local Server | AMD EPYC 7502 | 128 GB | NVIDIA Quadro RTX A6000 |
| Cloud Server | NVIDIA GH200 Grace Hopper Superchip | 463.9 GB | NVIDIA GH200 Grace Hopper Superchip |

**Appendix C. Results of minimizing the 10 benchmark test functions**

As detailed in Section 4, the GPU-based numerical method is employed to minimize the 10 benchmark test functions with scalable dimensions defined in SI Appendix A using four computational resources. The configurations of the four computational resources are provided in SI Appendix B. The results of the Levy function with different dimensions are summarized in Table 2. The results of the other nine benchmark functions with different dimensions are provided below, from Table S2 to Table S10. In these tables, "NA" represents that the corresponding computational resource is not used to minimize the benchmark function with the specific dimension because of an expected long execution time.

**Table S2.** Results of the Ackley function with different dimensions

| Dimension | Number of Iterations | Number of Output Regions | Computation Time | | | |
|---|---|---|---|---|---|---|
| | | | Laptop | Workstation | Local Server | Cloud Server |
| *n* = 50 | 50 | 1 | 383 s | 82 s | 180 s | 11 s |
| *n* = 100 | 100 | 1 | 1,508 s (0.42 h) | 321 s (0.09 h) | 701 s (0.19 h) | 29 s (0.01 h) |
| *n* = 500 | 500 | 1 | 38,127 s (10.59 h) | 8,062 s (2.24 h) | 17,604 s (4.89 h) | 565 s (0.16 h) |
| *n* = 1,000 | 1,000 | 1 | NA | 33,025 s (9.17 h) | 72,625 s (20.17 h) | 2,288 s (0.64 h) |

**Table S3.** Results of the Belegundu function with different dimensions

| Dimension | Number of Iterations | Number of Output Regions | Computation Time | | | |
|---|---|---|---|---|---|---|
| | | | Laptop | Workstation | Local Server | Cloud Server |
| *n* = 50 | 45 | 1 | 197 s | 42 s | 85 s | 6 s |
| *n* = 100 | 90 | 1 | 804 s (0.22 h) | 171 s (0.05 h) | 330 s (0.09 h) | 15 s |
| *n* = 500 | 450 | 1 | 24,728 s (6.87 h) | 5,276 s (1.47 h) | 8,919 s (2.48 h) | 426 s (0.12 h) |
| *n* = 1,000 | 900 | 1 | NA | 23,780 s (6.61 h) | 37,683 s (10.47 h) | 1,896 s (0.53 h) |

**Table S4.** Results of the Breiman function with different dimensions

| Dimension | Number of Iterations | Number of Output Regions | Computation Time | | | |
|---|---|---|---|---|---|---|
| | | | Laptop | Workstation | Local Server | Cloud Server |
| *n* = 50 | 40 | 1 | 203 s | 43 s | 135 s | 5 s |
| *n* = 100 | 80 | 1 | 792 s (0.22 h) | 166 s (0.05 h) | 531 s (0.15 h) | 13 s |
| *n* = 500 | 400 | 1 | 19,960 s (5.54 h) | 4,208 s (1.17 h) | 13,365 s (3.71 h) | 301 s (0.08 h) |
| *n* = 1,000 | 800 | 1 | NA | 17,041 s (4.73 h) | 53,621 s (14.89 h) | 1,218 s (0.34 h) |

**Table S5.** Results of the Fu function with different dimensions

| Dimension | Number of Iterations | Number of Output Regions | Computation Time | | | |
|---|---|---|---|---|---|---|
| | | | Laptop | Workstation | Local Server | Cloud Server |
| *n* = 50 | 45 | 1 | 737 s | 140 s | 373 s | 13 s |
| *n* = 100 | 90 | 1 | 2,636 s (0.73 h) | 555 s (0.15 h) | 1,472 s (0.41 h) | 41 s (0.01 h) |
| *n* = 500 | 450 | 1 | 65,867 s (18.30 h) | 13,844 s (3.85 h) | 36,688 s (10.19 h) | 892 s (0.25 h) |
| *n* = 1,000 | 900 | 1 | NA | 55,774 s (15.49 h) | 148,104 s (41.14 h) | 3,538 s (0.98 h) |

**Table S6.** Results of the Griewank function with different dimensions

| Dimension | Number of Iterations | Number of Output Regions | Computation Time | | | |
|---|---|---|---|---|---|---|
| | | | Laptop | Workstation | Local Server | Cloud Server |
| *n* = 50 | 55 | 1 | 416 s | 88 s | 194 s | 9 s |
| *n* = 100 | 110 | 1 | 1,653 s (0.46 h) | 348 s (0.10 h) | 763 s (0.21 h) | 28 s (0.01 h) |
| *n* = 500 | 550 | 1 | 43,256 s (12.02 h) | 9,077 s (2.52 h) | 19,958 s (5.54 h) | 605 s (0.17 h) |
| *n* = 1,000 | 1,100 | 1 | NA | 37,060 s (10.29 h) | 82,044 s (22.79 h) | 2,523 s (0.70 h) |

**Table S7.** Results of the Rastrigin function with different dimensions

| Dimension | Number of Iterations | Number of Output Regions | Computation Time | | | |
|---|---|---|---|---|---|---|
| | | | Laptop | Workstation | Local Server | Cloud Server |
| *n* = 50 | 45 | 1 | 276 s | 59 s | 179 s | 6 s |
| *n* = 100 | 90 | 1 | 1,086 s (0.30 h) | 229 s (0.06 h) | 703 s (0.20 h) | 18 s (0.01 h) |
| *n* = 500 | 450 | 1 | 27,178 s (7.55 h) | 5,706 s (1.59 h) | 17,678 s (4.91 h) | 400 s (0.11 h) |
| *n* = 1,000 | 900 | 1 | NA | 23,214 s (6.45 h) | 71,765 s (19.93 h) | 1,727 s (0.48 h) |

**Table S8.** Results of the Salomon function with different dimensions

| Dimension | Number of Iterations | Number of Output Regions | Computation Time | | | |
|---|---|---|---|---|---|---|
| | | | Laptop | Workstation | Local Server | Cloud Server |
| *n* = 50 | 55 | 1 | 157 s | 31 s | 70 s | 5 s |
| *n* = 100 | 110 | 1 | 641 s (0.18 h) | 120 s (0.03 h) | 265 s (0.07 h) | 12 s |
| *n* = 500 | 550 | 1 | 20,722 s (5.76 h) | 3,255 s (0.90 h) | 7,221 s (2.01 h) | 328 s (0.09 h) |
| *n* = 1,000 | 1,100 | 1 | NA | 13,606 s (3.78 h) | 31,409 s (8.72 h) | 1,435 s (0.40 h) |

**Table S9.** Results of the Styblinski function with different dimensions

| Dimension | Number of Iterations | Number of Output Regions | Computation Time | | | |
|---|---|---|---|---|---|---|
| | | | Laptop | Workstation | Local Server | Cloud Server |
| *n* = 50 | 45 | 1 | 268 s | 60 s | 138 s | 7 s |
| *n* = 100 | 90 | 1 | 1,055 s (0.29 h) | 233 s (0.06 h) | 572 s (0.16 h) | 21 s (0.01 h) |
| *n* = 500 | 450 | 1 | 26,741 s (7.43 h) | 7,090 s (1.97 h) | 15,210 s (4.23 h) | 419 s (0.12 h) |
| *n* = 1,000 | 900 | 1 | NA | 24,218 s (6.73 h) | 64,924 s (18.03 h) | 1,791 s (0.50 h) |

**Table S10.** Results of the Zabinsky function with different dimensions

| Dimension | Number of Iterations | Number of Output Regions | Computation Time | | | |
|---|---|---|---|---|---|---|
| | | | Laptop | Workstation | Local Server | Cloud Server |
| *n* = 50 | 40 | 1 | 399 s | 85 s | 189 s | 9 s |
| *n* = 100 | 80 | 1 | 1,578 s (0.44 h) | 331 s (0.09 h) | 734 s (0.20 h) | 26 s (0.01 h) |
| *n* = 500 | 400 | 1 | 41,658 s (11.57 h) | 8,808 s (2.45 h) | 19,366 s (5.38 h) | 619 s (0.17 h) |
| *n* = 1,000 | 800 | 1 | NA | 36,699 s (10.19 h) | 80,920 s (22.48 h) | 2,635 s (0.73 h) |

## Appendix D. Results of minimizing the Rosenbrock function with scalable dimensions

The Rosenbrock function with scalable dimensions has the form (15, 16):

$$f(\boldsymbol{x}) = \sum_{i=1}^{n-1} \left[ 100\left(x_i^2 - x_{i+1}\right)^2 + (x_i - 1)^2 \right], \tag{S22}$$

$$-30 \leq x_i \leq 30, \tag{S23}$$

where $i \in \{1, 2, \ldots, n\}$. The global minimum of the Rosenbrock function is located at $\boldsymbol{x}^* = [1, 1, \ldots, 1]$, and $f(\boldsymbol{x}^*) = 0$.

As detailed in Section 4, the GPU-based numerical method is employed to minimize the Rosenbrock function defined by Eq. (S22) and Eq. (S23) with 50, 100, 500, and 1,000 dimensions ($n$ = 50, 100, 500, and 1,000 in the function) using four different computational resources. The configurations of the four computational resources are provided in SI Appendix B. The results of this computational experiment are summarized in Table S11 below. In Table S11, "NA" represents that the corresponding computational resource is not used to minimize the Rosenbrock function with the specific dimension because of an expected long execution time. As shown in Table S11, the GPU-based numerical method outputs multiple regions when minimizing the Rosenbrock function with scalable dimensions, where each output region satisfies the region width tolerance of $1 \times 10^{-4}$. These output regions are interconnected and could be combined as one region that is guaranteed to include the global minimum of the Rosenbrock function.

**Table S11.** Results of the Rosenbrock function with different dimensions

| Dimension | Number of Iterations | Number of Output Regions | Computation Time | | | |
|---|---|---|---|---|---|---|
| | | | Laptop | Workstation | Local Server | Cloud Server |
| $n$ = 50 | 209 | 9 | 1,456 s | 309 s | 691 s | 26 s |
| $n$ = 100 | 595 | 18 | 8,259 s (2.29 h) | 1,743 s (0.48 h) | 3,908 s (1.09 h) | 118 s (0.03 h) |
| $n$ = 500 | 14168 | 99 | NA | NA | 467,761 s (129.93 h) | 12,577 s (3.49 h) |
| $n$ = 1,000 | 80548 | 267 | NA | NA | NA | 144,656 s (40.18 h) |

## Appendix E. Results of minimizing the Ackley function using different methods

As detailed in Section 4, seven popular CPU-based numerical optimization methods are employed to minimize the Ackley function with 100 dimensions using a laptop and a local server. The Ackley function is defined in SI Appendix A. The configurations of the laptop and local server are provided in SI Appendix B. The results of this computational experiment using the laptop are summarized in Table S12. The results using the local server are provided below in Table S13. Double-precision floating-point format (FP64 format) is employed for CPU and GPU computing in these methods. The minimum function value derived by each method is rounded to four significant digits to save space. The function value at the known global minimum is $f(\boldsymbol{x}^*) = 0$. All seven popular CPU-based numerical optimization methods fail to find the global minimum of the function, while the GPU-based numerical method (proposed method) successfully encloses the global minimum of the function in a single run.

**Table S12.** Results of the Ackley function with 100 dimensions ($n$ = 100) using the laptop

| Method | Number of Runs | Minimum Function Value | Total Computation Time |
|---|---|---|---|
| Basin-hopping | 1,000 | 6.472 | 6,899 s |
| BFGS | 50,000 | 7.091 | 3,020 s |
| Differential evolution | 100 | 0.7854 | 8,559 s |
| DIRECT | 1 | 3.326 | < 1 s |
| Interior-point | 50,000 | 5.807 | 2,205 s |
| Nelder-Mead | 10,000 | 7.489 | 6,967 s |
| Trust region | 50,000 | 4.403 | 3,528 s |
| **Proposed method** | **1** | **[$-5.684\times10^{-14}$, $3.469\times10^{-7}$]** | **1,508 s** |

**Table S13.** Results of the Ackley function with 100 dimensions ($n$ = 100) using the local server

| Method | Number of Runs | Minimum Function Value | Total Computation Time |
|---|---|---|---|
| Basin-hopping | 1,000 | 6.829 | 5,506 s |
| BFGS | 50,000 | 7.020 | 2,621 s |
| Differential evolution | 100 | 0.8985 | 12,008 s |
| DIRECT | 1 | 3.326 | < 1 s |
| Interior-point | 50,000 | 5.880 | 2,860 s |
| Nelder-Mead | 10,000 | 6.990 | 11,379 s |
| Trust region | 50,000 | 4.534 | 4,522 s |
| **Proposed method** | **1** | **[$-5.684\times10^{-14}$, $3.469\times10^{-7}$]** | **701 s** |

## Appendix F. Results of minimizing the discontinuous test function

A discontinuous test function with scalable dimensions is built from the Levy function (8) and the Dirac delta function (17):

$$f(\boldsymbol{x}) = \begin{cases} -1, & \text{if } 1.5 - \dfrac{\varepsilon}{2} \le x_i \le 1.5 - \dfrac{\varepsilon}{2}, \\ L(\boldsymbol{x}), & \text{otherwise,} \end{cases} \tag{S24}$$

where $i \in \{1, 2, \ldots, n\}$, $\varepsilon = 1\times10^{-6}$, and $L(\boldsymbol{x})$ represents the Levy function:

$$L(\boldsymbol{x}) = \frac{\pi}{n}\left\{10\sin^2(\pi y_1) + (y_n - 1)^2 + \sum_{i=1}^{n-1}\left[(y_i - 1)^2\left(1 + 10\sin^2(\pi y_{i+1})\right)\right]\right\}, \tag{S25}$$

$$y_i = 1 + 0.25(x_i - 1), \tag{S26}$$

$$-10 \le x_i \le 10. \tag{S27}$$

The global minimum of the discontinuous test function is introduced by the Dirac delta function with $f(\boldsymbol{x}^*) = -1$. The one-dimensional discontinuous test function ($n$ = 1 in the function) is illustrated in Figure S1.

As detailed in Section 4, the GPU-based numerical method is employed to minimize the discontinuous test function defined by Eqs. (S24) - (S27) with 50, 100, 500, and 1,000 dimensions ($n$ = 50, 100, 500, and 1,000 in the function) using four different computational resources. The configurations of the four computational resources are provided in SI Appendix B. The results of this computational experiment are summarized in Table S14 below. In Table S14, “NA” represents that the corresponding computational resource is not used to minimize the discontinuous test function with the specific dimension because of an expected long execution time.

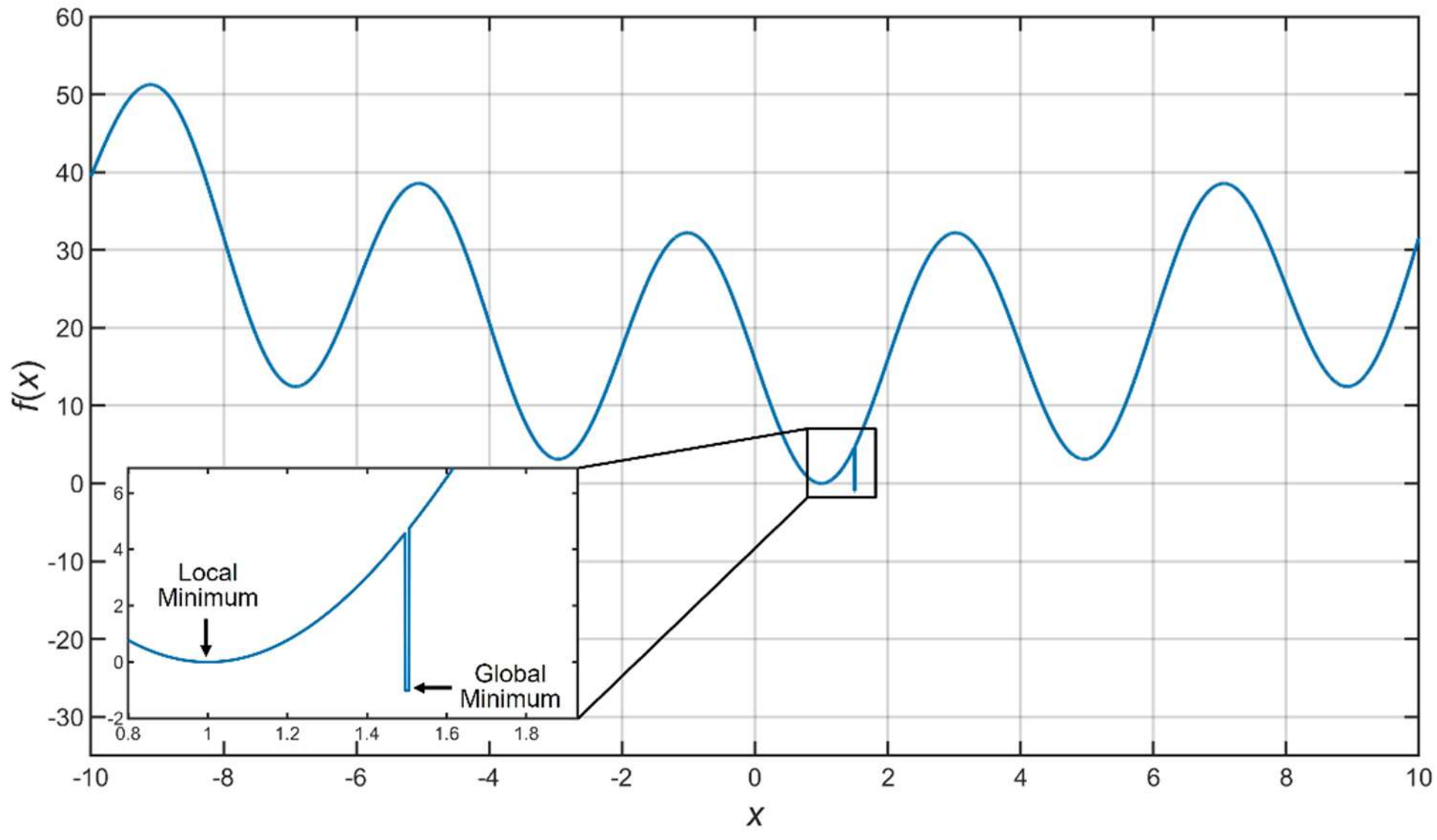

**Figure S1.** The discontinuous test function with one dimension ($n$ = 1)

**Table S14.** Results of the discontinuous test function with different dimensions

| Dimension | Number of Iterations | Number of Output Regions | Computation Time | | | |
|---|---|---|---|---|---|---|
| | | | Laptop | Workstation | Local Server | Cloud Server |
| $n$ = 50 | 45 | 1 | 638 s | 366 s | 832 s | 35 s |
| $n$ = 100 | 90 | 1 | 2,512 s (0.70 h) | 1,454 s (0.40 h) | 3,289 s (0.91 h) | 114 s (0.03 h) |
| $n$ = 500 | 450 | 1 | 62,445 s (17.35 h) | 34,950 s (9.71 h) | 82,451 s (22.90 h) | 2,527 s (0.70 h) |
| $n$ = 1,000 | 900 | 1 | NA | 135,184 s (37.55 h) | NA | 9,928 s (2.76 h) |

**SI References**